\documentclass[11pt]{amsart}
\usepackage{amsfonts,amsmath,amssymb}
\usepackage[all]{xy}
\usepackage{hyperref}
\begin{document}

% 30 July 2020

\newcommand{\ie}{\textit{i}.\textit{e}.\,}
\newcommand{\eg}{\textit{e}.\textit{g}.\,}
\newcommand{\cf}{{\textit{cf}.\;}}
\parindent=0pt 
\parskip=6pt

\newcommand{\la}{{\langle}}
\newcommand{\ra}{{\rangle}}
\newcommand{\alg}{{\rm alg}}
\newcommand{\A}{{\mathcal A}}
\newcommand{\B}{{\bf B}}
\newcommand{\BFK}{{\rm BFK}}
\newcommand{\bk}{{\boldsymbol{\kappa}}}
\newcommand{\bo}{{\boldsymbol{\omega}}}
\newcommand{\bt}{{\boldsymbold{\vartheta}}}
\newcommand{\bbar}{{\bar{\omega}}}
\newcommand{\ba}{{\bf a}}
\newcommand{\bc}{{\bf c}}
\newcommand{\bi}{{\bf i}}
\newcommand{\bI}{{\bf I}}
\newcommand{\bq}{{\bf k}}
\newcommand{\bm}{{\bf m}}
\newcommand{\bn}{{\bf n}}
\newcommand{\bu}{{\bf u}}
\newcommand{\bv}{{\bf v}}
\newcommand{\bY}{{\bf Y}}
\newcommand{\bh}{{h_{M\U}}}
\newcommand{\bx}{{h_{M\xi}}}
\newcommand{\br}{{\rm BR}}
\newcommand{\bDelta}{{\boldsymbol{\Delta}}}
\newcommand{\cZ}{{\mathcal Z}}
\newcommand{\C}{{\mathbb C}}
\newcommand{\cHd}{{{\mathcal H}^{dif}}}
\newcommand{\cHi}{{{\mathcal H}^{inv}}}
\newcommand{\M}{{\mathbb M}}
\newcommand{\N}{{\mathbb N}}
\newcommand{\Q}{{\mathbb Q}}
\newcommand{\R}{{\mathbb R}}
\newcommand{\mS}{{\sf S}}
\newcommand{\bS}{{\bf S}}
\newcommand{\Spec}{{\rm{sp} \;}}
\newcommand{\T}{{\mathbb T}}
\newcommand{\Z}{{\mathbb Z}}
\newcommand{\Tor}{{\rm Tor}}
\newcommand{\Hom}{{\rm Hom}}
\newcommand{\Cay}{{\rm Cay}}
\newcommand{\Lie}{{\rm Lie \;}}
\newcommand{\Sp}{{\rm Sp}}
\newcommand{\U}{{\rm U}}
\newcommand{\op}{{\rm op}}
\newcommand{\dR}{{\rm dR}}
\newcommand{\res}{{\rm res}}
\newcommand{\cQ}{{\sf Q}}
\newcommand{\cN}{{\sf N}}
\newcommand{\sN}{{\bf N}}
\newcommand{\SU}{{\rm SU}}
\newcommand{\sslash}{{/\!\!/}}

\title{ Renormalization groupoids in algebraic topology}

\author[Jack Morava]{Jack Morava}

\address{Department of Mathematics, The Johns Hopkins University,
Baltimore, Maryland 21218}

\email{jack@math.jhu.edu}

% \subjclass{Primary ... ; Secondary ...}                                     
  
% \keywords{...}                                                              
  
% \thanks{...}                                                                
  
% \dedicatory{...}                                                            
  
% \date{August 2020}

\begin{abstract}{Continuing work begin in \cite{19,27}, we interpret the 
Hurewicz homomorphism for Baker and Richter's noncommutative complex cobordism 
spectrum $M\xi$ in terms of characteristic numbers (indexed by quasi-symmetric 
functions) for complex-oriented quasitoric manifolds, and show that 
automorphisms or cohomology operations on this representation are defined 
by the `renormalization' Hopf algebra $\sN_*$ of formal diffeomorphisms at
the origin of the noncommutative line, previously considered (over $\Q$) 
in quantum electrodynamics \cite{3}.\medskip 

\noindent The resulting structure can be presented in purely algebraic terms, 
as a groupoid scheme over $\Z$ defined by a coaction of $\sN_*$ on the ring
$\cN_*$ of noncommutative symmetric functions. We sketch some applications 
to symplectic toric manifolds, combinatorics of simplicial spheres, and 
statistical mechanics.}\end{abstract} \bigskip

\maketitle \bigskip

{\bf Introduction}

In \cite{1,2} Andrew Baker and Birgit Richter defined a remarkable 
noncommutative analog $M\xi$ of the complex cobordism spectrum $M\U$.
The homology $H_*(M\U,\Z)$ is a free polynomial algebra with one generator
in each even degree, while $H_*(M\xi,\Z)$ is a free associative algebra
with one generator in each even degree. The homology of $M\U$ can be
identified with the algebra of `characteristic numbers' of complex-oriented
manifolds (defined for example by the Chern-Weil theory of integrals of
curvature forms), while $H_*(M\xi,\Z)$ can be identified with invariants
defined by quasi-symmetric functions of the line bundles associated to 
complex-oriented quasi-toric manifolds. 

In the late 1960s, work of Landweber and Novikov led to an analysis of 
the cohomology operations in complex cobordism in terms of what has come 
to be understood as a `Hopf algebroid' structure
\[
M\U_*M\U \cong M\U_* \otimes \bS_* \;,
\]
where $\bS_*$ is the commutative but noncocommutative Hopf algebra 
representing the group of formal diffeomorphisms of the line at 
the origin, known previously (over $\Q$) to combinatorialists as the 
Faa di Bruno algebra. An easy corollary identifies 
\[
H_*(M\U \wedge M\U,\Z) \cong H_*(M\U,\Z) \otimes \bS_*
\]
in terms of that algebra, and the ring of characteristic numbers. 

It is fair to say that the noncommutative ring $M\xi_*M\xi$ defining 
the cohomology operations in $M\xi_*$ is not yet well-understood, but
Baker and Richter show that it injects into $H_*(M\xi \wedge M\xi,\Z)$,
and that this injection becomes an isomorphism over $\Q$. The purpose
of this note is to identify 
\[
H_*(M\xi \wedge M\xi.\Z) \cong H_*(M\xi,\Z) \otimes \sN_*
\]
in terms of a remarkable generalization $\sN_*$ of $\bS_*$, neither
commutative nor cocommutative, already known (over $\Q$) to physicists
\cite{3,4} as the `renormalization' Hopf algebra of formal diffeomorphisms
at the origin of the noncommutative line. 

Since the work of Quillen it has become clear that it is useful to think
of complex cobordism in terms of the moduli object for one-dimensional
formal group laws. As a first step toward a deeper understanding of 
$M\xi$, we suggest an interpretation of the algebra $H_*(\C P^\infty 
\wedge M\xi,\Z)$ as functions on a noncommutative `grand canonical 
ensemble' of characteristic numbers of symplectic quasitoric manifolds,
with morphisms defined by noncommutative renormalization. \bigskip

{\sc \S I \; Some algebra}\bigskip

[This section reviews the basic properties of the algebra of symmetric
functions (used in the theory of characteristic classes) following 
\cite{21}, and of the Hopf algebra representing the group of formal 
diffeomorphisms of the line at the origin used in cobordism theory \cite{31}. 
It then summarizes some work of Brouder, Frabetti, Krattenthaler, and Schmitt 
\cite{3,4} on a (neither commutative nor cocommutative) Hopf algebra 
representing formal diffeomorphisms of the {\bf non}commutative line, 
used in renormalization theory.] \bigskip

{\bf \S 1 The unique self-dual irreducible positive graded Hopf algebra} 

{\bf 1.1} Let $\Z[x_*] = \Z_{k \geq 1}[x_k]$ be the graded polynomial 
algebra on generators $x_k$ of degree $|x_k| = 2$; then the subalgebra 
\[
\Z[e_*] = \Z_{k \geq 1}[e_k] := \mS_* \subset \Z[x_*]
\]
of elementary symmetric `functions' \cite{21}( I \S 2.7 p 22) is generated by 
elements $e_k$ of degree $2k$, defined by the formal series
\[
e(T) = \prod_{k \geq 1} (1 + x_k T) = \sum_{j \geq 0} e_j T^j \in \mS_*[[T]]\;.
\]
Similarly, 
\[
h(T) = e(-T)^{-1} = \prod_{k \geq 1}(1 - x_k T)^{-1} = \sum_{j \geq 0} h_j T^j
\in \mS_*[[T]] \; 
\]
defines the `complete' symmetric functions of degree $|h_k| = 2k$; thus 
$e_0 = h_0 = 1, \; e_1 = h_1 \dots$. If $I = 1^{i_1}2^{i_2} \cdots r^{i_r}$ is 
an (unordered) partition of $|I| = \sum_{1 \leq k \leq r} k i_k$ into $r(I) = 
r$ parts, let 
\[
e_I := \prod_{1 \leq k \leq r} e_k^{i_k} \in \mS_{2|I|} \;,
\]
and similarly for $h_I$, etc. Elements of the form $e_I$ (or $h_I$ or $p_I$, 
see \cite{21}(I \S 6) and further below) provide bases for $\mS_*$.

The coproduct $\Delta_\mS : \mS_* \to \mS_* \otimes_\Z \mS_*$ defined by 
\[
\Delta e(T) = e(T) \otimes e(T) 
\]
(undecorated $\otimes$ signifies $\otimes_\Z$) and the antipode 
\[
\chi_\mS(e_n) = \sum_{|I| = n} (-1)^{r(I)} e_I
\]
make $\mS_*$ into a binomial Hopf algebra. The group-valued Witt functor
\[
A \mapsto (\Spec \mS)(A) = \Hom_\alg(\mS,A)
\]
on commutative rings (ignoring the grading)\begin{footnote}{sp 
is the canonical contravariant functor from a category to its 
opposite}\end{footnote} sends a ring homomorphism $\alpha : \mS \to A$ 
to the (multiplicatively invertible) power series
\[
\alpha(e)(T) = \sum_{k \geq 0} \alpha(e_k)T^k \in (1 + TA[[T]])^\times \;;
\]
the lost grading can be recovered from an action of the multiplicative
group, as suggested below.\bigskip

{\bf 1.2} P Hall's positive definite inner product \cite{21}(\S 1.4 p 63, 1.5 
ex 25 p 91) on $\mS_*$ defines a canonical isomorphism with its dual Hopf 
algebra, rendering some applications to topology confusing. The (primitive)  
power sums
\[
p(-T) := \sum_{k \geq 1} p_k T^{k-1} = e(T)^{-1} \cdot e'(T) \in \mS_*[[T]]
 \]
define an orthogonal basis $p_I$ for $\mS_* \otimes \Q$, while the monomial 
symmetric functions
\[
m_I := \sum_{\sigma \in \Sigma_r}  x^{i_1}_{\sigma(1)} \cdots 
x^{i_r}_{\sigma(r)} \in \mS_{2|I|}
\]
are dual to the $h_I$. The cohomology of the classifying space for stable 
complex vector bundles (see II 1.2) can be identified with $\mS_*$, with 
the Chern class $c_k \in H^{2k}(B\U,\Z)$ corresponding to $e_k$. 

The operation which sends a stable vector bundle to its dual sends $e_k$ 
to $(-1)^ke_k$, while the Whitney sum inverse $V \mapsto -V$ interchanges 
$e_k$ and $(-)^kh_k$, defining a kind of adjoint involution $\omega : V 
\mapsto -V^*$. Hall's duality $h_I \leftrightarrow m_I$ is a Hopf algebra 
involution on $\mS_*$. \bigskip

{\bf \S 2. Formal diffeomorphisms of the line at 0}

{\bf 2.1} When $A$ is a commutative algebra, the set 
\[
(\Spec \bS)(A) = \{t(T) = \sum_{k \geq 0} t_k T^{k+1} \in A[[T]]\:|\: t_k \in
A \}
\] 
is a noncommutative monoid under composition $t',t \mapsto (t' \circ t)(T) = 
t'(t(T))$ of formal power series. If $t_0$ is invertible, an elementary 
induction shows that $t$ has a unique two-sided compositional inverse 
$t^{-1}$ with $(t^{-1} \circ t)(T) = (t \circ t^{-1})(T) = T$. It follows 
that the localized polynomial algebra
\[
\bS_* := t_0^{-1} \Z_{k \geq 0}[t_k] := \Z[t_0^{\pm 1}] \otimes \tilde{\bS}_*
\] 
has a noncocommutative coproduct which can be expressed as
\[
(\Delta_\bS t)(T) = (t \otimes 1)((1 \otimes t)(T)) \in (\bS_* \otimes 
\bS_*)[[T]] \;;
\]
the composition inverse defines an antipode $\chi_\bS$ on $\bS_*$ given 
explicitly by Lagrange's reversion formula, defining a commutative but 
noncocommutative Hopf algebra. Enlarging $\tilde{\bS}_*$ to $\bS_*$ encodes 
the grading as a coaction of the Hopf algebra of the multiplicative 
groupscheme (with coproduct $\Delta t_0 = t_0 \otimes _0$). This can be 
suppressed by taking $t_0 = 1$.

The noncommutative dual algebra $\bS^*$ appears in algebraic topology as the 
analog of the Steenrod algebra for complex cobordism, where it is known as the 
Landweber-Novikov algebra. Quillen's theorem interprets $\bS_*$ as the Hopf 
algebra representing the group of coordinate transformations of 
one-dimensional formal group laws. 

The noncommutative groupscheme $\Spec \bS = \mathbb{G}_m \times \Spec 
\tilde{\bS}$ (\ie of formal diffeomorphisms at the origin of the line 
over $\Spec \Z$) is closely related to the combinatorialists' Faa 
di Bruno algebra: or, more precisely, to the analogous group of formal 
diffeomorphisms of the line over $\Q$, but now parametrized by 
Taylor-MacLaurin coordinates $t^{(n)}(0) = n! t_{n-1}$. 

{\bf 2.2} There is a right action 
\[
\Spec \mS \times_{\Spec \Z} \Spec \bS \to \Spec \mS
\]
of the noncommutative groupscheme $\Spec \bS$ on the commutative groupscheme 
$\Spec \mS$, defined by composition 
\[
e(T),t(T) \mapsto (e \circ t)(T)
\]
of power series, \ie by $\psi_\mS(e(T)) = (e \otimes 1)(t(1 \otimes T))$,
representing a coaction 
\[
\psi_\mS : \mS_* \to \mS_* \otimes \bS_* 
\]
making $\mS_*$ into a Hopf $\bS_*$ - algebra comodule. This structure can be
interpreted as defining a `right unit' $\eta_R = \psi_\mS$ and a 
`left unit' $\eta_L = {\bf 1}_\mS \otimes 1$ for a (split) Hopf algebroid
\cite{31} 
\[
\mS_\bullet \bS := (\mS_*, \mS_* \otimes \bS_*) : \mS_* \Rightarrow \mS_* 
\otimes \bS_*
\]
with coproduct 
\[
(\mS_* \otimes \bS_*) \to \mS_* \otimes (\bS_* \otimes \bS_*) \cong 
(\mS_* \otimes \bS_*) \otimes_{\mS_*} (\mS_* \otimes \bS_*) \;. 
\] 
These morphisms extend to define a cosimplicial (Amitsur) algebra
\[
\A^\bullet_{\mS_*}(\mS_* \otimes \bS_*) : 0  \to \mS_* \Rightarrow 
\mS_* \otimes \bS_* \Rrightarrow \mS_* \otimes (\bS_* \otimes \bS_*) \dots
\;,
\]
a cobar construction or codescent complex \cite{16}(\S 5)  
\[
\A^\bullet_C(B) : 0 \to C \Rightarrow B \Rrightarrow B \otimes_C B \dots
\]
such that $D \otimes_C \A^\bullet_C(B) \cong \A^\bullet_D (D \otimes_C B)$, \eg
\[
\A^\bullet_{\mS_*}(\mS_* \otimes \bS_*) \cong \mS_* \otimes \A^\bullet_\Z
(\bS_*) \;.
\]
Such constructions underly (for example) the Adams-Novikov spectral sequences
of homotopy theory. In the dual language of schemes
\[
\Spec \A^\bullet_{\mS_*} (\mS_* \otimes \bS_*) \simeq [\Spec \mS \sslash 
\Spec \bS]
\]
represents a kind of pre-stack: the (untopologized) groupoid-valued functor on 
commutative algebras defined by the action of coordinate transformations on
the symmetric functions.\bigskip

{\bf \S 3. Noncommutative and quasisymmetric functions}

{\bf 3.1} The free associative graded cocommutative (but noncommutative) 
binomial Hopf algebra 
\[
\cN_* = \Z_{k \geq 1} \langle Z_k \rangle = \Z \langle Z_* \rangle
\]
of noncommutative symmetric functions\begin{footnote}{As in 2.1, the grading
on $\cN_*$ can be encoded by adjoining a central unit $Z_0$ of degree 
zero}\end{footnote} (\ie with coproduct $\Delta Z(T) = Z(T) \otimes Z(T), 
\; Z(T) = \sum_{k \geq 0} Z_k T^{k+1}$) is dual to a commutative but not 
cocommutative algebra 
\[
\cQ_* \subset \Z[x_*] 
\]
of {\bf quasi}symmetric functions: an {\bf ordered} partition 
\[
\bI = i_1 + \cdots + i_r 
\]
of the integer $|I|$ into $r$ {\bf nonempty} parts (there are $2^{|I|-1}$ 
such things) defines a basis element \cite{1},\cite{15}(\S 4) 
\[
m_\bI = \sum_{0 < n_1 < \cdots < n_r} \prod_{1 \leq j \leq r} x_{n_j}^{i_j} 
\in \cQ_{2|I|} \subset\Z[x_*] 
\]
dual to the basis element $Z_\bI = Z_{i_1} \cdots Z_{i_r} \in \cN_{2|I|}$. 
We have 
\[
\chi_\cN(Z_n) = \sum_{|I| = n} (-1)^{r(I)} Z_\bI
\]
for the antipode of $\cN$. According to Ditters' conjecture, $\cQ_*$ is 
a free commutative algebra; over $\Q$ it is polynomial, generated 
by elements of degree $2n$ indexed by Lyndon words of degree $n$, with
multiplicative structure defined by a certain shuffle product on the 
basis $\{m_\bI\}$. The quasisymmetric functions have
interesting relations with Koszul duality \cite{24}(\S 3) and number 
theory \cite{6}(\S 4.4), but these will not be pursued here. The surjection
\[
\cN_* \ni Z_k \mapsto e_k \in \mS_* 
\]
of Hopf algebras defined by abelianization dualizes to a monomorphism 
$\mS_* \to \cQ_*$.\bigskip

{\bf 3.2} Although composition of power series over noncommutative rings is 
{\bf not} associative in general \cite{2}, there is a remarkable generalization
of the Hopf algebra of formal diffeomorphisms of the line in the noncommutative
context. Brouder, Frabetti, and Krattenthaler \cite{3}(Thm 2.4) define a Hopf 
algebra $\sN_*$ with coproduct 
\[
(\Delta_\sN Z)(T) = \res_{U=0} \; Z(U) \otimes_\Z (U - Z(T))^{-1} \in (\sN 
\otimes_\Z \sN)_*[[T]]
\]
on generators of the underlying algebra $\cN_*$ via a formal version 
\[
\res_{U=0} \; (\prod_{i \in \Z} a_i U^i) := a_{-1}
\]
of Cauchy's residue\begin{footnote}{Here and throughout, variables such 
as $T,U,\dots$ will always be central, of cohomological degree 
two.}\end{footnote} which elegantly simplifies proof of coassociativity: 
expanding the right-hand side, we have 
\[
\res_{U=0} \; U^{-1} Z(U) \otimes_\Z (1 - U^{-1} Z(T))^{-1} 
= \sum_{k \geq 0} \res_{U=0} \; U^{-k-1} Z(U) \otimes Z(T)^k 
\]
\[
= \sum_{k \geq 1} Z_{k-1} \otimes_\Z Z(T)^k = (Z \otimes 1)((1 \otimes Z)(T))
\;, 
\]
\ie
\[
\Delta_\sN Z_k = \sum_{|I|+j=k} Z_j \otimes Z_\bI \;.
\]
Their Theorem 2.14 provides an explicit formula for the 
antipode $\chi_\sN$ of the resulting Hopf algebra. Similarly, their \S 3 
shows that the coproduct, regarded as an algebra homomorphism
\[
\Delta_\sN : \cN_* \to \cN_* \otimes \sN_* \;,
\]
is compatible with the binomial coalgebra structure on $\cN_*$, making it
a Hopf algebra comodule over $\sN_*$ \cite{18}(III \S 7). As in \S 2.2, this 
can be reformulated as the assertion that
\[
\cN_\bullet \sN := (\cN_*,\cN_* \otimes \sN_*) : \cN_* \Rightarrow 
\cN_* \otimes \sN_*
\]
is a (noncommutative, split) Hopf algebroid over $\Z$, with $\mS_\bullet \bS$
as its abelianization. \bigskip

{\bf \S 4 A digression, on renormalization Hopf algebras}

{\bf 4.0} I owe Michiel Hazewinkel thanks for drawing attention to the Hopf 
algebra $\sN_*$ and its applications in the formalization of quite classical 
quantum electrodynamics. It is a part of a wider literature (largely over 
$\Q$) expressed in terms of combinatorial Hopf algebras of trees, graphs, 
posets, and related structures, and to place it appropriately in that enormous 
field is not practical here; but its many possible generalizations deserve at 
least some mention. This section collects a few key ideas from work of Brouder 
and Schmitt on a general class of `renormalization bialgebras' $T(T(\B)^+)$ 
constructed functorially as cotensor algebras on a graded bialgebra $\B$. When 
$\B = \Z$ is the trivial bialgebra, 
\[
T(\B)^+ = \{ \otimes^n e \:|\: n \geq 1\} 
\]
is the augmentation ideal of the (co)tensor algebra $T(\B) = \Z[e]$, 
identifying $Z_n$ with $(\otimes^n e) \in T^{2n}(T(\B))$. \bigskip

{\bf 4.1} A certain partially ordered monoid $C$ of `compositions' -- ordered 
partitions as in \S 3.1, but with $\emptyset$ allowed as identity 
element -- plays a basic role in \cite{4}: concatenation, \eg
\[
(1 + 1) * (2 + 3) = (1 + 1 + 2 + 3)
\]
is its multiplication operation. The iterated cotensor algebra $T(T(\B)^+)$
is shown to be naturally $C$-graded, and Brouder and Schmitt associate to 
$\rho, \sigma \in C$, certain restriction and contraction coalgebra 
morphisms
\[
T^\sigma(T(\B)^+) \to T^{\sigma|\rho}(T(\B)^+) \;, 
\]
\[
T^\sigma(T(\B)^+) \to T^{\sigma/\rho}(T(\B)^+)
\]
of a sort familiar (for example) in the study of Hopf algebras of trees
and graphs. On generators $u \in T^\tau(\B))$ with 
\[
\Delta u = \sum u'_i \otimes u''_i \;,
\]
they define a new, `renormalization' coproduct
\[
\bDelta u = \sum_{\sigma \leq \tau} u'_i|\sigma \otimes u''_i/\sigma \;, 
\]
which then extends to define a new bialgebra structure \cite{3}(\S 2.3 Th 1)
on $T(T(\B)^+)$. \bigskip

{\bf 4.2} Constructions of this sort are ubiquitous in work on combinatorial
Hopf algebras, but vary in details, which are omitted here. If $\varepsilon$
is the counit of $\B$, then elements of the form $(x) - \varepsilon(x)$
generate a two-sided ideal $(J) \subset T(T(\B)^+)$, such that the quotient
bialgebra $T(T(\B)^+)/(J)$ is a Hopf algebra \cite{3}(\S 5.1). This ideal is
trivial for $\B = \Z$, and so recovers (in different notation) the
construction of $\sN_*$. \bigskip

{\sc \S II Some topology} \bigskip

{\bf \S 1 The Hurewicz homomorphism}

[This section reviews the theory of characteristic classes and numbers for
complex cobordism in the language of \S I. Its point is that the 
noncommutative cobordism spectrum $M\xi$ is not yet well-understood, but 
the associated theory of characteristic numbers is surprisingly accessible.]
\bigskip

{\bf 1.1} A topological group $G$, regarded as a topological category and 
hence as a simplicial space, has a canonical classifying space $BG$ for 
$G$-bundles as its topological realization. For example, the circle group 
$\T = \U(1)$ has B$\T \simeq \C P^\infty$, with polynomial cohomology 
\[
H^*(\C P^\infty,\Z) \cong \Z[c], \; |c| = 2
\]
generated by the Chern class of a complex line bundle. $\T$ being commutative, 
the product map $\T \times \T \to \T$ is a group homomorphism, making 
$H^*(B\T,\Z)$ a primitively generated Hopf algebra. Its dual 
\[
H_*(B\T,\Z) = \Z_{k \geq 1}[b_{(k)}]
\]
is a divided power algebra, with generators satisfying
\[
b_{(n)} \cdot b_{(m)} = (n,m) \: b_{(n+m)} \;.
\]

A theorem of IM James identifies the stable homotopy type of the
loopspace $\Omega \Sigma X$ (of the reduced suspension of a simply-connected
space) with that of the free topological monoid
\[
\Omega \Sigma X \sim \bigvee_{k \geq 1} X^{\wedge k}
\]
generated by the pointed space $X$. [The free {\bf commutative} topological
monoid on $X$ is the infinite symmetric product ${\rm SP}^\infty X$ of Thom and
Dold.] For example, a level one projective representation of the loop group
$L\SU(2))$ defines a morphism 
\[
\Omega \SU(2) \cong \Omega \Sigma S^2 \to \C P^\infty
\]
of loop spaces, with an induced morphism
\[
H_*(\Omega S^3,\Z) = \Z[\beta] \to \Z[b_{(*)}] = H_*(\C P^\infty,\Z)
\]
of Hopf algebras sending $\beta \in \pi_2(\Omega S^3)$ to $b_{(1)}$. 
[This notation is meant to suggest Boltzmann's thermodynamic beta; see 
1.3.2 below.]

Similarly, the homology of $\C P^\infty$ is torsion-free, so the K\"unneth 
theorem implies that 
\[
H_*(\Omega \Sigma \C P^\infty,\Z) \cong \Z \langle Z_* \rangle \cong \cN_* 
\]
is isomorphic to the cofree cotensor coalgebra on $\tilde{H}_*(\C P^\infty,
\Z)$, with generator $b_{(k)}$ corresponding to $Z_k$, defining a basis 
$Z_\bI$ as in I 3.1 and \cite{1}.\bigskip

{\bf 1.2} Block diagonal or Whitney sum composition
\[
\U(n) \times \U(m) \to \U(n+m)
\]
defines a homotopy commutative- and -associative law on the monoid
\[
\coprod_{n \geq 0} B\U(n) \;,
\]
whose group completion
\[
\Omega B (\coprod_{n \geq 0} B\U(n)) \simeq \Z \times B\U \supset 0 \times B\U
\]
has a classifying space for stable complex vector bundles as identity 
component. The maximal toruses $\T^n \subset \U(n)$ then define Borel's Hopf
algebra isomorphiam
\[
H^*(B\U,\Z) \cong \mS_* \;.
\]
\cite{5}, while the map $B\U(1) \to B\U$ defines a polynomial basis
\[
H_{2k}(B\U(1),\Z) \ni b_{(k)} \mapsto b_k \in H_{2k}(B\U,\Z) 
\]
with $b(T) = \sum_{i \geq 0}b_kT^{k+1}$ satisfying $\Delta b(T) = b(T) \otimes
b(T)$. The associated primitives $\check{p}_k$ (as in I 1.2) then satisfy 
$c_i \cdot \check{p}_k = \delta_{i,k}$.  

The universal complex vector bundle $\xi_n \to B\U(n)$ defines the Thom
spectrum 
\[
M\U := \{ S^n \wedge \xi^+_m = S^n \wedge M\U(n) \to M\U(n+m) = \xi^+_{n+m} \}
\]
representing the complex cobordism ring
\[
M\U_* = \pi_* M\U = \lim_{k \to \infty} \pi_{*+k} M\U(k)
\]
of compact complex-oriented manifolds. Hurewicz's ring homomorphism
\[
\pi_* M\U \to H_*(M\U,\Z) \cong H_*(B\U,\Z) \cong \mS_*
\]
takes the homotopy groups of the composition 
\[
M\U \simeq M\U \wedge S^0 \to M\U \wedge H\Z 
\]
of morphisms of ring spectra defined by the identity map
\[
[S^0 \to H\Z] = 1 \in \tilde{H}^0(S^0,\Z) \;.
\]
(where $S^0$ is the sphere spectrum and $H\Z$ is the integral
Eilenberg - Mac Lane spectrum). This extends to a natural characteristic 
number homomorphism
\[
h_{M\U} : M\U^*(X) \to H^*(X,\mS_*)
\]
of multiplicative cohomology functors, familiar (over $\R$) from Chern-Weil 
theory.\bigskip

{\bf 1.3.1} A cobordism class $\{X\} \in M\U_{2k}$ can be interpreted as the 
equivalence class of a compact smooth submanifold $X \subset \R^{2k+N}, \; 
N \gg 0$, with complex normal bundle classified by 
\[
\nu_X : X \to B\U(N) \to B\U \;.
\]
The Thom - Pontryagin collapse construction
\[
S^* \to \R^*/(\R^* - \nu(X)) = \nu(X)^+ \to \xi^+_*
\]
defines a map
\[
\xymatrix{ 
H_*(S^*) \ar[r]^-{\nu_{X*}} & H_*(\xi^+_*) \ar[r] & H_*(M\U) \cong H_*(B\U) 
\cong \mS_* }
\]
sending the generator of $H_*(S^*)$ to the Hurewicz image of $\{X\}$. A dual 
construction sends $\{X\}$ to the system of integrals of the tangential Chern 
classes $c_I(\tau_X)$ with dim $X = 2|I|$, regarded as a linear functional on
$\mS_*$. Hall duality equates this with the linear functional 
\[
h_{M\U}\{X\} : m_I \mapsto (-1)^{|I|}m_I(\nu_X) \cap [X] \in \mS_*
\]
defined by the monomial characteristic classes of the normal bundle. 

{\bf 1.3.2} Because $M\U_* \to H_*(M\U,\Z)$ is a homomorphism of torsion-free 
rings, it is convenient to work with their rationalizations. Quillen's theorem
identifies the completed Hopf algebra $M\U^*\C P^\infty$ with Lazard's 
universal one-dimensional formal group law, and over $\Q$ such group laws are 
classified by their logarithms. Mi\v{s}\v{c}enko's theorem 
\[
\log_{M\U}(T) = \sum_{k \geq 1} \{\C P_{k-1}\} \frac{T^k}{k} \in M\U^*_\Q[[T]]
\]
identifies this logarithm; on the other hand
\[
\bh : M\U^*(\C P^\infty) \to H^*(\C P^\infty,H_*(B\U,\Z)) \cong H^*(\C 
P^\infty,\mS_*)
\]
expresses the formal group law in terms of the exponential
\[
\bh (X +_{M\U} Y) = b(b^{-1}(X) + b^{-1}(Y)) \in \mS_*[[X,Y]]
\]
defined $b(T)$ as above. Note that 
\[
h_{M\U}\{\C P_n\} = (n+1) \: \chi_\mS(b_n) \;,
\]
see \cite{18}(28.5.012).

In \cite{32} Ravenel and Wilson show that the canonical inclusions 
\[
\beta_k = \{\C P_k \subset \C P^\infty\} \in M\U_{2k}(\C P^\infty)
\]
generate the complex bordism of $\C P^\infty$, modulo the relation
\[
\beta(X +_{M\U} Y) = \beta(X) \cdot \beta(Y)
\]
(with $\beta(T) = \sum_{k \geq 0}\beta_k T^{k+1}$). Working rationally, it 
follows that
\[
\log (\beta \circ b)(T) = b_{(1)}T 
\]
and hence (recalling that $b_{(1)}$ is the image of $\beta$ as in 1.4)  
\[
\beta(T) = \exp (\beta \log_{M\U}(T)) \in (\mS^\Q_*[\beta])[[T]] \;.
\]
As suggested by Friedrich and McKay \cite{11}(Prop 4.2) this resembles 
formally the canonical partition function in statistical mechanics;
see further in \S III. 

{\bf 1.4} As in I 2.1, the moduli object $[\Spec M\U_* \sslash \Spec S_*]$ for 
one-dimensional group laws can thus be identified with the stack defined by 
the now-classical Hopf algebroid
\[
M\U_\bullet M\U : M\U_* \Rightarrow M\U_*M\U \;;
\]
the Hurewicz or characteristic number map 
\[
h_{M\U} : [\Spec H_*M\U \sslash \Spec \bS_*] \subset [\Spec M\U_* \sslash 
\Spec \bS_*]  
\]
then becomes the inclusion of the stratum of formal group laws of 
additive type. It pulls back to an isomorphism over $\Spec \Q \subset 
\Spec \Z$. \bigskip

{\bf \S 2 The Baker-Richter spectrum} $M\xi$

{\bf 2.1} Following 1.2, group completion of the monoidal map
\[
\coprod_{n \geq 0} B\U(1)^{\times n} \to \coprod_{n \geq 0} B\U(n)
\]
defines a morphism $\br : \Omega \Sigma B\U(1) \to \Z \times B\U$ of $A_\infty$
spaces. Pulling back the stable universal bundle $\xi$ over $0 \times B\U$
defines Baker and Richter's $A_\infty$ spectrum
\[
M\xi := \{ S^n \wedge \br^*\xi^+_m \to \br^*\xi^+_{n+m} \}
\]
together with an abelianization morphism $M\xi \to M\U$. A remarkable
theorem \cite{1}(\S 7) shows that localization at a prime splits $M\xi$ as 
a wedge of suspensions of the Brown-Peterson spectrum, even though $M\xi$ is 
not an $M\U$-module spectrum; but it {\bf is} complex-orientable, and possesses
natural Thom isomorphisms.

We will however make no use of this local structure here; our focus is the 
ring-spectrum $M\xi \wedge H\Z$, \ie on the integral homology and
cohomology
\[
H_*(M\xi,\Z) \cong \cN_*, \; H^*(M\xi,\Z) \cong \cQ_*
\]
and their relation to $M\xi$ through its Hurewicz homomorphism. Baker
and Richter show that 
\[
M\xi_* = \pi_*M\xi \to H_*(M\xi,\Z)
\]
is injective, and an isomorphism after rationalization. Following 1.3 and 
\cite{19}, we can interpret the Hurewicz map as a noncommutative 
characteristic number homomorphism
\[
\bx : M\xi^*(Y) \to H^*(Y,H_*(M\xi,\Z)) \cong H^*(Y,\cN_*) \cong \Hom(\cQ_*,
H^*(Y,\Z)) \;.
\]
Under Quillen's conventions \cite{30} we can  assume that $Y$ is a smooth 
manifold, and interpret an element $\{X\} \in M\xi^{2i}(Y)$ as the cobordism 
class of a map $[\nu : X \to Y]$ between manifolds, with stable normal bundle 
\[
v_X \cong \oplus L_i : X \to \bigvee_{n \geq 0} B\U(1)_+^{\times}
\]
presented as a direct sum of stable complex line bundles. Then $\bx \{X\}$ 
is represented by the linear functional
\[
m_\bI \mapsto (-1)^i \nu_!(m_\bI(c(L_*))) \in H^{2i}(Y,\Z)
\]
on $\cQ_{2i}$ (where $\nu_!$ is the covariant pushforward cohomology 
homomorphism induced by the complex-oriented map $\nu$).  

{\bf Remark} The right unit 
\[
\eta_R : \sN_* \to \cN_* \otimes \cN_* \cong \Hom(\cQ_*,\cN_*)
\]
of I \S 2.3 can be defined on generators by
\[
\eta_R(Z_k)(m_\bI) = \sum_{i+j=k} m_\bI(Z_i) \cdot Z_j \;, 
\]
agreeing with the coproduct formula in I 3.2. \bigskip

{\bf 2.2} The Euler class $x_\xi \in M\xi^2 \C P^\infty$ provides a
(noncentral \cite{2}) coordinate for kind of noncommutative one-dimensional 
formal group structure on $M\xi^*\C P^\infty$. Its Hurewicz image 
\[
h_{M\xi}(x_\xi) = \sum_{k \geq 0}Z_kT^{k+1} = Z(T) \in H_*(M\xi,\Z)[[T]] 
\subset \cN_*[[T]]
\]
defines a (completed) Hopf coproduct 
\[
\bDelta_{\C P^\infty}(T) = \sum_{k \geq 0} Z_k \cdot (\chi_\sN Z)(T \otimes 1)
+ (\chi_\sN Z)(1 \otimes T))^{k+1} \in \cN_*[[T \otimes 1,1 \otimes T]]\;,
\]
\cf Zassenhaus's noncommutative binomial theorem. 

The Thom isomorphism $M\xi_*M\xi \cong M\xi_* \Omega \Sigma \C P^\infty$
(together with a K\"unneth argument for torsion-free spaces) then defines
a cosimplicial noncommutative algebra 
\[
\A^\bullet_{M\xi_*}(M\xi_*M\xi) : 0 \to M\xi_* \Rightarrow M\xi_*M\xi 
\Rrightarrow M\xi_*M\xi \otimes_{M\xi_*} M\xi_*M\xi \dots
\]
as in I 1.2. In view of the remark above, the Hurewicz morphism
\[
M\xi \wedge M\xi \to (M\xi \wedge H\Z) \wedge_{H\Z} (M\xi \wedge H\Z) \simeq
M\xi \wedge M\xi \wedge H\Z
\]
induces an injective homomorphism
\[
\A^\bullet_{M\xi_*}(M\xi_*M\xi) \to \A^\bullet_{H_*M\xi}(H_*M\xi \otimes 
H_*M\xi) \cong \A^\bullet_{\cN_*}(\cN_* \otimes \sN_*)\;,
\]
of cosimplicial algebras defined by homotopy groups: dually, a morphism 
\[
[\Spec \cN \sslash \Spec \sN] \to \Spec \A^\bullet_{M\xi_*}(M\xi_*M\xi)  
\]
of pre-stacks, which becomes an isomorphism over $\Q$. Thinking of 
commutative objects as a subclass of noncommutative ones defines an
abelianization morphism 
\[
[\Spec \mS \sslash \Spec \bS] \to [\Spec \cN \sslash \Spec \sN] \;.
\]
It is tempting to imagine the descent object $\A^\bullet_{M\xi}(M\xi \wedge
M\xi)$ as an approximation to a noncommutative analog of the sphere spectrum. 
\bigskip

{\sc \S III Examples and remarks} 

{\bf 1} A $2n$-dimensional prequantized toric manifold \cite{9,12,14,27} is an 
omnioriented complex-oriented manifold (\ie with a preferred decomposition 
of its stable tangent bundle as a sum of complex lines) defined by an effective
action of a torus $\T^n$, together with an equivariant complex line bundle
$L$ with connection, such that the curvature form $\omega(\nabla_L)$ is
symplectic; it is roughly analogous to a compactified system of harmonic 
oscillators. 

Such a manifold defines an element of $M\xi_{2n}\C P^\infty$
and thus has characteristic numbers in $H_*(\C P^\infty, \cN_*) \cong 
\cN_*[b_{(*)}]$. Similarly toric varieties, considered as orbifolds 
\cite{19,20}, have characteristic numbers in $\cN_*\otimes \Q$, while 
Hamiltonian toric varieties have similar invariants in $M\xi_*\C P^\infty 
\otimes \R \cong M\xi_* \otimes \R[\beta]$. Projective toric varieties 
\cite{13}
provide examples; in particular, chemical reaction networks \cite{8}
with $k$ nodes and deficiency zero \cite{27}(\S 1.3.2) have characteristic 
invariant $\sum_{\bI=k} Z_\bI$. 

{\bf 2} The equivariant cohomology of the moment-angle complex associated
to a simplicial (\eg Gorenstein \cite{27}(\S 1.2.2)) complex can be 
identified with its Stanley-Reisner face ring, which thus has 
(quasi-symmetric) characteristic classes and numbers. This defines a
homomorphism \cite{19} from the (noncommutative) ring of simplicial spheres,
with join as composition, to $\sN_*$; the boundary of the $k$-simplex yields
the example above. 

{\bf 3} Work in the theory of free, \ie noncommutative, probability \cite{28},
\cite{17}(Th 3.1),\cite{29} suggests that
\[
(\chi_\sN Z)(-T) \in \sN_*[[T]]
\]
generalizes the cumulant generating function \cite{10,22} of classical 
statistical mechanics, along the lines suggested in II 1.3.2 above. 

\newpage

\bibliographystyle{amsplain}

\begin{thebibliography}{99}

\bibitem[1]{1} A Baker, B Richter, Quasisymmetric functions from a topological
point of view, Math. Scand. 103 (2008) 208 -– 242,
\url{https://arxiv.org/abs/math/0605743}

\bibitem[2]{2} -----, -----, Some properties of the Thom spectrum over loop
suspension of complex projective space, in {\bf An alpine expedition through   
algebraic topology} 1 -– 12, Contemp. Math. 617, AMS 2014,
\url{https://arxiv.org/abs/1207.4947}

\bibitem[3]{3} C Brouder, A Frabetti, C Krattenthaler, Non-commutative Hopf
algebra of formal diffeomorphisms, Adv. Math. 200 (2006) 479 –- 524,
\url{https://arxiv.org/abs/math/0406117}

\bibitem[4]{4} -----, W Schmitt, Renormalization as a functor on bialgebras,
J Pure Appl. Algebra 209 (2007) 477 -– 495,
\url{https://arxiv.org/abs/hep-th/0210097}

\bibitem[5]{5} H Cartan, D\'emonstration homologique des th\'eor\'em\`es de
p\'eriodicit\'e de Bott, II, Séminaire Henri Cartan, tome 12 (1959 -- 1960)
exp.17, \url{http://www.numdam.org/volume/SHC_1959-1960__12/}

\bibitem[6]{6} P Cartier, A primer of Hopf algebras, in {\bf Frontiers in 
number theory, physics, and geometry II}  537 -– 615, Springer (2007),
\url{http://preprints.ihes.fr/2006/M/M-06-40.pdf}

\bibitem[7]{7} P Conner, E Floyd, {\bf The relation of cobordism to 
$K$-theories}, Springer LNM 28 (1966)

\bibitem[8]{8} G Craciun, A Dickenstein, A Shiu, B Sturmfels, Toric
dynamical systems, J. Symbolic Comput. 44 (2009) 1551 –- 1565,
\url{https://arxiv.org/abs/0708.3431}

\bibitem[9]{9} M Davis, T Januszkiewicz, Convex polytopes, Coxeter orbifolds
and torus actions, Duke Math. J. 62 (1991) 417 –- 451

\bibitem[10]{10} R Friedrich, J McKay, Free probability and complex cobordism,
C R Math Acad Sci So R Can 33 (2011) 116 –- 122

\bibitem[11]{11} ------, -----, Formal groups, Witt vectors, and free
probability, \url{https://arxiv.org/abs/1204.6522}

\bibitem[12]{12} VL Ginzburg, Calculation of contact and symplectic
cobordism groups, Topology 31 (1992) 767 -– 773

\bibitem[13]{13} V Guillemin, Kaehler structure on toric varieties,
J Diff Geo 40 (1994) 285 - 309


\bibitem[14]{14} V Guillemin, VL Ginzburg, Y Karshon, {\bf Moment maps,
cobordisms, and Hamiltonian group actions}, AMS Mathematical Surveys and
Monographs 98 (2002)

\bibitem[15]{15} M Hazewinkel, Symmetric functions, noncommutative symmetric
functions, and quasisymmetric functions, Acta Appl. Math. 75 (2003) 55 -- 83,
\url{https://arxiv.org/abs/0410468}

\bibitem[16]{16} K Hess, A general framework for homotopic descent and
codescent, \url{https://arxiv.org/abs/1001.1556}

\bibitem[17]{17} M Josuat-Verg\'es, F Menous, JC Novelli, JY Thibon, Free
cumulants, Schroder trees, and operads, Adv in App Math 88 (2017) 92 -- 119,
\url{https://arxiv.org/abs/1604.04759}

\bibitem[18]{18} C Kassel, {\bf Quantum groups}, Springer Graduate Texts 155
(1995)1995.

\bibitem[19]{19} N Kitchloo, J Morava, The Baker-Richter spectrum as cobordism
of quasitoric manifolds, \url{https://arxiv.org/abs/1201.3127}

\bibitem[20]{20} E Lerman, S Tolman, Hamiltonian torus actions on symplectic
orbifolds and toric varieties, Trans. AMS 349 (1997) 4201 -– 4230,
\url{https://arxiv.org/abs/dg-ga/9511008}

\bibitem[21]{21} I Macdonald, {\bf Symmetric functions and Hall polynomials},
2nd ed, Oxford Classic Texts in the Physical Sciences (2008)

\bibitem[22]{22} M Marcolli, Motivic information,
\url{https://arxiv.org/abs/1712.08703}

\bibitem[23]{23} A Mathew, Climbing Mount Bourbaki,
\url{https://amathew.wordpress.com/2012/05/28/}

\bibitem[24]{24} J Morava, Homotopy-theoretically enriched categories of
noncommutative motives, Res. Math. Sci. 2 (2015), Art. 8, 16 pp.,
\url{https://arxiv.org/abs/1402.3693}

\bibitem[25]{25} -----, On formal groups and geometric quantization,
\url{https://arxiv.org/abs/1905.06181}

\bibitem[27]{27} -----, Topological invariants of some chemical reaction
networks, \url{https://arxiv.org/abs/1910.12609}

\bibitem[28]{28}  A Nica, R Speicher, {\bf Lectures on the combinatorics of
free probability}, LMS  Lecture Note Series 335, Cambridge (2006)

\bibitem[29]{29} JC Novelli, JY Thibon, Hopf algebras and dendriform
structures arising from parking functions, Fund. Math. 193 (2007) 189 – 241,
\url{https://arxiv.org/abs/math/0511200}

\bibitem[30]{30} D Quillen, Elementary proofs of some results of cobordism
theory using Steenrod operations, Adv. in Math. 7 (1971) 29 -– 56

\bibitem[31]{31} D Ravenel, {\bf Complex cobordism and stable homotopy groups  
of spheres}, Pure and Applied Mathematics 121, Academic Press (1986)

\bibitem[32]{32} -----, W S Wilson, The Hopf ring for complex cobordism,
J Pure Appl Algebra 9 (1976/77) 241 -– 280

\end{thebibliography}

\end{document}